\documentclass[11pt]{amsart}
\usepackage[linktocpage]{hyperref}
\usepackage{amssymb, paralist, xspace, graphicx, url, amscd, euscript, mathrsfs,stmaryrd,epic,eepic,color}
\usepackage[all]{xy}
\usepackage{amsthm}
\usepackage{enumerate}
\usepackage{lipsum}
\usepackage[ampersand]{easylist}
\SelectTips{cm}{}
\usepackage{multirow}

\numberwithin{equation}{section}

1
\numberwithin{subsection}{section}

\allowdisplaybreaks[1]

\newtheorem*{namedtheorem}{\theoremname}
\newcommand{\theoremname}{testing}

\newtheorem{theorem}[subsection]{Theorem}
\newtheorem{thm}[subsection]{Theorem}
\newtheorem{proposition}[subsection]{Proposition}

\newtheorem{proposition-definition}[subsection]
{Proposition-Definition}
\newtheorem{corollary}[subsection]{Corollary}
\newtheorem{lemma}[subsection]{Lemma}
\theoremstyle{definition}
\newtheorem{definition}[subsection]{Definition}
\newtheorem{question}[subsection]{Question}

\newtheorem{remark}{Remark}

\theoremstyle{remark}

\newcommand\cC{\mathcal{C}}
\newcommand\cD{\mathcal{D}}

\newcommand\cF{\mathcal{F}}

\newcommand\cK{\mathcal{K}}

\newcommand\cM{\mathcal{M}}

\newcommand\cO{\mathcal{O}}
\newcommand\cP{\mathcal{P}}

\newcommand\cR{\mathcal{R}}

\newcommand\cT{\mathcal{T}}
\newcommand\cU{\mathcal{U}}
\newcommand\cV{\mathcal{V}}
\newcommand\cW{\mathcal{W}}

\newcommand\CC{\mathbb{C}}

\newcommand\NN{\mathbb{N}}

\newcommand\PP{\mathbb{P}}
\newcommand\QQ{\mathbb{Q}}

\newcommand\ZZ{\mathbb{Z}}

\newcommand\rank{{\rm rank}}
\newcommand\Pic{{\rm Pic}}

\newcommand\IGr{{\rm IGr}}
\newcommand\LGr{{\rm LGr}}
\newcommand{\q}{/\!\!/}

\theoremstyle{plain}
\theoremstyle{definition}

\begin{document}

\title{Picard groups on moduli of K3 surfaces with Mukai models}
\author{Francois Greer,~~Zhiyuan Li,~~Zhiyu Tian}

\address{Department of Mathematics\\
building 380\\
Stanford, CA 94305\\
U.S.A.}

\address{Department of Mathematics\\
building 380\\
Stanford, CA 94305\\
U.S.A.}

\address{Department of Mathematics\\California Institute of Technology\\ Pasadena, CA,
91125\\
U.S.A.}
\email{fgreer@math.stanford.edu,~~zli2@stanford.edu,~~tian@caltech.edu}
\begin{abstract}
We discuss the Picard group of the moduli space $\cK_g$ of quasi-polarized K3 surfaces of genus $g\leq 12$ and $g\neq 11$. In this range, $\cK_g$ is unirational, and a general element in $\cK_g$ is a complete intersection with respect to a vector bundle on a homogenous space, by the work of Mukai.  In this paper, we find generators for the Picard group $\Pic_\QQ(\cK_g)$ using Noether-Lefschetz theory.  This verifies the Noether-Lefschetz conjecture on the moduli of K3 surfaces in these cases.
\end{abstract}
\maketitle

\section{Introduction}
It is well-known that the moduli space $\cM_g$  of smooth projective curves of genus $g$ is a quasi-projective variety with Picard number one for $g\geq 3$ (cf.~\cite{Ha83}). The Picard group $\Pic_\QQ(\cM_g)$ with rational coefficients is generated by the first Chern class of the Hodge bundle. In the moduli theory of higher dimensional varieties, the primitively polarized K3 surface of genus $g$ can be viewed as a two dimensional analogue of genus $g$ smooth projective curve, and it is natural for us to study the Picard group of its moduli space $\cK_g^\circ$.

Unlike in the curve case, the rank of $\Pic(\cK^\circ_g)$ is no longer constant;  it has been shown by O'Grady \cite{OG86} that $\rank( \Pic(\cK^\circ_g)) $ increases to infinity as $g$ is increasing (see also \cite{MP07} $\S 7$). Besides the Hodge line bundle, there are many other natural divisors on $\cK_g^\circ$ coming from geometry. Actually, the Noether-Lefschetz locus of $\cK^\circ_g$ parametrizing K3 surfaces with Picard number greater than one is a union of countably many irreducible divisors by Hodge theory. We call them the {\it Noether-Lefschetz (NL) divisors} on $\cK_g^\circ$. In this paper, we study the Picard group on $\cK^\circ_g$ for $g\leq 12$ ($g\neq 11$) and find its generators in terms of NL-divisors on $\cK_g^\circ$.

For conveniency, instead of working on $\cK^\circ_g$, we will study the moduli space $\cK_{g}$ of primitively quasi-polarized K3 surfaces of genus $g$, which is more natural from a Hodge theoretic point of view.  It is known that  $\cK_g$ is a locally Hermitian symmetric variety, by the Torelli theorem, and the complement $\cK_g\backslash \cK_g^\circ$ is a divisor parametrizing $K3$ surfaces containing a $(-2)$-exceptional curve.  In this setting, we define the Noether-Lefschetz (NL) divisors on $\cK_g$ as follows: given $d,n\in\ZZ$, let $D_{d,n}^{g}\subset \cK_{g}$ be the locus of those K3 surfaces $S\in \cK_g$ whose Picard lattice $\Pic(S)$ contains a primitive rank two sublattice
\begin{equation}\label{lattice}\begin{array}{c|c|c} & L & \beta \\\hline L & 2g-2 & d \\\hline \beta & d & n\end{array}
\end{equation}
where $L$ is the primitive quasi-polarization of $S$ and $\beta\in \Pic(S)$. Each NL divisor $D_{d,n}^g$ is irreducible by \cite{O86}.  In \cite{MP07}, Maulik and Pandharipande conjectured that the divisors $\{D_{d,n}^g\}$ span the  group $\Pic_\QQ(\cK_g)$ with rational coefficients. Our main result  is: 
\begin{theorem}
The Picard group of  $\cK_{g}$ with rational coefficients is spanned by NL divisors for $g\leq 10$ and $g=12$. Moreover, the basis of $\Pic_\QQ(\cK_{g})$ of $6\leq g\leq 10,$ is given by:
\begin{itemize}
\item $g=6$, $\{D^{6}_{0,-2}, ~D^{6}_{5,2}, ~D^{6}_{k,0},~k=1,\ldots 4\}$. 
\item $g=7$, $\{D^{7}_{0,-2}, ~D^{7}_{5,2}, ~D^{7}_{6,2},~D^{7}_{k,0},~k=1,\ldots 4\}$.
\item $g=8$, $\{D^{8}_{0,-2}, ~D^{8}_{6,2}, ~D^{8}_{7,2}~D^{8}_{k,0},~k=1,\ldots 4\}$.
\item $g=9$, $\{D^{9}_{0,-2}, ~D^{9}_{6,2}, ~D^{9}_{7,2}~D^{9}_{k,0},~k=1,\ldots 5\}$.
\item $g=10$, $\{D^{10}_{0,-2}, ~D^{10}_{7,2},~D^{10}_{8,2},~ D^{10}_{9,4}, ~D^{10}_{k,0},~k=1,\ldots 5\}$.
\end{itemize}
For $g=12$, the group $\Pic_\QQ(\cK_{12})$ is generated by 
\begin{equation}\label{generator}
\{D^{12}_{0,-2}, ~D^{12}_{7,2},~D^{12}_{8,2},~D^{12}_{9,2},~ D^{12}_{10,4},~ D^{12}_{11,4}, ~D^{12}_{k,0},~k=1,\ldots 6\}.
\end{equation}
\end{theorem}
\begin{remark}
The rank of $\Pic(\cK_{12})$ is $11$, so there is a linear relation between the generators in \eqref{generator}. See Remark \ref{relation} for more details.
\end{remark}
When $2\leq g \leq 5$, a general K3 surface in $\cK_g$ is a double cover of $\PP^2$ ($g=2$) or a complete intersection in $\PP^{g}$, and the assertion can be found in \cite{O86}\cite{Sh80}\cite{Sh81}\cite{LT13}. In these cases, the proof relies on the explicit construction of the moduli space of the corresponding complete intersections via geometric invariant theory (GIT). If $g$ is greater than $5$, the general K3 surface of genus $g$ is no longer a complete intersection in $\PP^g$, but it can be interpreted as a complete intersection with respect to certain vector bundles on homogenous spaces, for $g$ in the range of our theorem, so there ought to be a similar construction.

\noindent{\bf Acknowledgments} We have benefited from conversations with Brendan Hassett and Jun Li. We are very grateful to Arie Peterson for his useful comments and providing us the relation of NL-divisors in Remark \ref{relation}. 

\section{Geometry of K3 surfaces}
In this section, we review Mukai's work on the projective models of general low genus quasi-polarized K3 surfaces. We give a precise characterization, in terms of the Picard lattice, of the (non-general) K3 surfaces lying outside the locus of Mukai models. We also include a few examples to illustrate the phenomenon.  Throughout this paper, we work over complex numbers.

\subsection{Preliminary} Let $(S,L)$ be a primitively quasi-polarized K3 surface of genus $g$, i.e. $L\in \Pic(S)$ is big and nef with $L^2=2g-2$. The middle cohomology $H^2(S,\ZZ)$ is an even unimodular lattice of signature $(3,19)$ under the intersection form $\left<,\right>$.  Let $\Lambda_{g}:=c_{1}(L)^{\perp}$ be the orthogonal complement, which is an even lattice of signature $(2,19)$.  The period domain $\cD_{g}$ associated to $\Lambda_{g}$ can be realized as a connected component of 
 $$\cD_{g}^{\pm}:=\{v\in \PP(\Lambda_{g}\otimes \CC)| \left<v,v \right>=0,-\left<v, \bar{v}\right> > 0\}.$$
The monodromy group 
$$\Gamma_{g}=\{g\in \text{Aut}(\Lambda_{g})^{+}|~g~\hbox{acts trivially on}~\Lambda^\vee_{g}/\Lambda_{g} \},$$
naturally acts on $\cD,$ where $\text{Aut}(\Lambda_{g})^{+}$ is the identity component of $\text{Aut}(\Lambda_{g})$.
According to the Global Torelli theorem for K3 surfaces, there is an isomorphism $$\cK_{g}\cong \cD_{g}/\Gamma_g$$ via the period map. Hence, $\cK_{g}$ is a locally Hermitian symmetric variety with only finite quotient singularities, and is thus $\QQ$-factorial. Each NL divisor $D^g_{d,n}$ can be considered as the quotient of a codimension one subdomain in $\cD$. 

The NL divisors we defined here are slightly different from the divisors discussed in \cite{MP07}.   Indeed, Maulik and Pandharipande define the NL divisors by specifying a curve class as below.  
Given the data $d,n$, let $\cC^g_{d,n}$ be the locus parametrizing the K3 surfaces containing a class $\beta$ with $\beta^2=n$ and $L\cdot\beta=d$. Then $\cC_{d,n}^g$ is a divisor and it is supported on a collection of NL divisors $D_{d',n'}^g$ satisfying $\Delta^g_{d',n'}=k^2\Delta^g_{d,n}$ for some $k\in \ZZ$, where $\Delta^g_{d,n}$ is the determinant of the lattice $\Lambda_{d,n}^g$.
It is not hard to see that the span of two sets $\{D^g_{d,n}\}$ and $\{\cC^g_{d,n}\}$ are the same, and we have $\cK_g\backslash\cK_g^\circ=\cC^g_{0,-2}$.

\subsection{Mukai models}
Let $(S,L)$ be a smooth quasi-polarized K3 surface of genus $g$. The linear system $|L|$ defines a map  $\psi_L$ from $S$ to $\PP^{g}$. The image of $\psi_L$ is called a {\it projective model} of $S$. 
If $\psi_L$ is birational to its image, the $\psi_L(S)$ is a degree $2g-2$ surface in $\PP^{g}$ with at worst rational double points. 
\begin{remark}\label{singular} Suppose $\psi_L$ is birational. If the K3 surface $S$ contains  a $(-2)$-exceptional curve, the morphism $\psi_L$ contracts this exceptional curve and the image $\psi_L(S)$ has a rational double point. Thus it is easy to see that $\psi_L(S)$ is singular only if $(S,L)\in \cK_g\backslash\cK_g^\circ$.
\end{remark}
As we mentioned in the introduction, when $6\leq g\leq 12$ ($g\neq 11$), general members in $\cK_{g}$ are no longer complete inspections in $\PP^{g}$, but can be interpreted as {\it complete intersections with respect to vector bundles} on homogenous spaces in the following sense.  Let  $E$ be a rank $r$ vector bundle on a smooth projective variety $X$ with local frame $\{e_1,\ldots, e_r\}$.  A global section $$s=\sum\limits_{i=1}^r f_ie_i\in H^0(X,E),~f_i\in \cO_X,$$  is non-degenerate at $x$ if $s(x)=0$ and $(f_1,\ldots, f_r)$ is a regular sequence. The section $s$ is non-degenerate if it is nondegenerate at every point of $(s)_0$, where  $(s)_0$ the zero locus of the section $s$. 

A subscheme $Y$ is a {\it complete intersection with respect to $E$} if $Y$ is the zero locus of a non-degenerate global section of $W$. In particular, one can view the complete intersection of $r$-hypersurfaces $H_1, H_2,\ldots, H_r$ in $X$ as the complete intersection with respect to the vector bundle $$E=\cO_X(H_1)\oplus \cO_X(H_2)\ldots \oplus \cO_X(H_r)$$ on $X$. 

Now we review Mukai's Brill-Noether theory on K3 surfaces, which essentially gives us the classification of the projective models of general K3 surfaces.
\begin{definition}
A polarized K3 surface $(S,L)$ of genus $g$ is Brill-Noether (BN) general if the inequality $h^0(M)h^0(N)<h^0(L)$ holds for any pair $(M,N)$ of non-trivial line bundles such that $M\otimes N\cong L$. 
\end{definition}

When $(S,L)$ is BN general of genus $g$, for any two integers $r,s$ with $rs=g$, Mukai \cite{Mu02} shows that there exists a rigid and stable vector bundle $E_{r}$ on $S$ of rank $r$ such that $\wedge^r E_{r}\cong L$.  The higher cohomology of $E_{r}$ vanishes and  $\dim H^0(S,E_{r})=r+s$.   Then there  is a  map
$$\Phi_{E_r}: S\rightarrow Gr(r, H^0(S, E_{r})^\vee). $$
\begin{remark}
The BN theory on K3 surfaces is an analogous to the BN theory on curves. Actually, if a smooth curve $C\in |L|$ is BN general, then the polarized K3 surface $(S,L)$ is BN general (cf.~\cite{Ha02}). 
\end{remark}

For $6\leq g \leq 12$ and $g \neq 11$, Mukai has shown that  image of $\Phi_{E_r}$ can be described as complete intersections with respect to a vector bundle on some homogenous subspace in $Gr(r, H^0(S,E_{r}))^\vee$. Here we just summarize Mukai's results in \cite{Mu02}:

\begin{thm}\label{mukaimodel}
A primitively quasi-polarized K3 surface $(S,L)$ of genus $6\leq g\leq 10$  or $g=12$ is Brill-Noether general if and only if it is birational to a complete intersection with respect to a vector bundle $E_{g}$ in a homogenous space $X_{g}$ via $\psi_L$, where  $X_{g}$ and the images $\psi_L(S)$ are listed below:
\begin{itemize}
\item$g=6$: $X_{6}=Gr(2,5)$ and $E_{6}=\cO_{X_{6}}(1)^{\oplus 3}\oplus \cO_{X_{6}}(2)$, $\psi_L(S)$ is a complete intersection of  a quadric  and  a codimension three linear section in $X_{6}$;

\item$g=7$:  $X_{7}$ is the isotropic Grassmannian  $\IGr(5,10)$ and $E_{7}=\cO_{X_{7}}^{\oplus 8}$, $\psi_L(S)$ is a codimension eight linear section of $X_{7}$;

\item$g=8$:  $X_{8}=Gr(2,6)$ and $E_{8}=\cO_{X_{8}}^{\oplus 6}$, $\psi_L(S)$ is a codimension six linear section of  $X_{8}$;

\item$g=9$:  $X_{9}$ is the Langrangian Grassmannian $\LGr(3,6)$ and $E_{9}=\cO_{X_{9}}^{\oplus 4}$, $\psi_L(S)$ is a codimension four linear section of  $X_{9}$;

\item$g=10$:  $X_{10}$ is the flag variety of dimension five associated with the adjoint representation of $G_2$  imbedded in $\PP^{13}$ and $E_{10}=\cO_{X_{10}}(1)^{\oplus 3}$;  $\psi_L(S)$ is a codimension three linear section of $X_{10}$;

\item $g=12$: let  $V$ be a seven dimensional vector space and $X_{12}=Gr(3,V)$. Then  $(S,L)\in \cK_{12}$ is BN general if and only if $\psi_L(S)$  is birational to  a hyperplane of  $Gr(3,V,N)\subseteq \PP^{13}$, where $N\subseteq \wedge^2 V^\vee$ is a non-degenerate three dimensional 
subspace and $Gr(3,V,N)\subseteq Gr(3,V)$ consists of  three dimensional subspaces  $U$ of $V$ such that the restriction of $N$ to $U\times U$ is zero. 
\end{itemize}
\end{thm}
Here, the non-degeneracy of the subspace $N\subseteq \wedge^2 V^\vee$ means that there is no decomposable vector in $N\wedge V^\vee\subseteq \wedge^3 V^\vee$. 
When $N$ is non-degenerate,  $Gr(3,V,N)$ is a smooth Fano threefold of index 1 and can be considered as  a complete intersection with respect to the vector bundle $\wedge^2 \cF^{\oplus 3}$ on $Gr(3,V)$, where $\cF$ is the dual of the universal subbundle on $Gr(3,V)$. 

\begin{remark}
For $g=6,8,9,10,12$, the values of $r,s$ chosen by Mukai  are 
\begin{center}
    \begin{tabular}{|l||l|l|l|l|l|}
    \hline
    $g$ & 6& 8 & 9 & 10 & 12 \\ \hline
    $rs$ & $2\cdot3$ &  $2\cdot 4$& $3\cdot3$ & $2\cdot5$ & $3\cdot 4$ \\ \hline
    \end{tabular}
\end{center}
In the case $g=7$, Mukai chose the rank 5 vector bundle $E_{12}$ with $\wedge^5E_{7}=L^{\otimes 2}$. The linear system $(E_{7}, H^0(S,E_7))$ embeds the K3 surface $S$ into a Grassmannian $Gr(5,10)$, whose image is contained in the isotropic Grassmannian $\IGr(5,10)$. 
\end{remark}

\begin{remark}The homogenous space $X_g$ is the quotient of a simply connected semisimple Lie group by a maximal parabolic subgroup. Here, we list some of the associated semisimple Lie groups of $X_{g}$ ($7\leq g\leq 10$), which will be used later:
\begin{itemize}[-]
  \item $g=7$: the spin group $Spin(10)$.
  \item $g=8$: the special linear group $PGL(6)$.
  \item $g=9$: the symplectic group $Sp(6)$.
  \item $g=10$: the exceptional group of type $G_2$.
\end{itemize}
\end{remark}

Next, observing that a smooth codimension three linear section of $Gr(2,5)$ is the unique (up to isomorphism) Fano threefold $\cF_5$ of degree five and index two, one can certainly consider the general K3 surface in $\cK_6$ as a quadric hypersurface in $\cF_5$.  It is known that the Fano threefold $\cF_5$ is a quasi-homogenous space with automorphism group $PSL(2)$.  The following Lemma  describes the locus of such K3 surfaces.
\begin{lemma}
A BN general K3 surface $(S,L)$ of genus $6$ is contained in the smooth Fano threefold $\cF_{5}$ if and only if $(S,L)$ is not contained in $D^{6}_{4,0}$. 
\end{lemma}
\begin{proof}
By Theorem \ref{mukaimodel}, we already know that a BN general K3 surface $(S,L)\in\cK_6$ is a complete intersection of a quadric hypersurface $Y\subseteq X_{6}$ and a codimension three linear section $\Sigma \subseteq X_{6}$. We know that $\Sigma$ is isomorphic to $\cF_5$ if $\Sigma $ is smooth. 

If $\Sigma$ is singular, we claim that $\Sigma$ must contain a quadric surface $Q$. Actually,
the codimension three linear sections in $X_{10}$ have been determined by Todd and Kimura (cf.~\cite{Todd},\cite{Ki82}), and it is straightforward to see the existence of the quadric surface. For instance, a codimension three linear section $\Sigma_0$ with an $A_1$ singularity (this is the generic case) is defined in $\PP^{6}$ as follows:
\begin{equation}
\begin{aligned}
z_0z_2+z_1z_4+z_3^2=0\\
z_0z_5+z_4^2+z_3z_6=0\\
z_1z_5-z_2z_4=0\\
z_1z_6-z_3z_4=0\\
z_2z_6-z_3z_5=0
\end{aligned}
\end{equation}
Then there exists a quadric cone $Q=\{z_0z_2+z_3^2=z_4=z_5=z_6=0\}$ in $\Sigma_0$.  The K3 surface $(S,L)$ is contained in $D^{6}_{4,0}$ since the intersection $Q\cap Y$ is an elliptic curve of degree four.
\end{proof}

\subsection{Non-BN general K3 surfaces}
 In this subsection, we classify all non-BN general K3 surfaces for $10\leq g \leq 12$ ($g\neq 11$) and interpret them as a union of NL divisors in $\cK_{g}$.   This is natural because non-BN general K3 surfaces $(S,L)$ must contain some special curve, and hence lie in some NL divisor $D^{g}_{d,n}$.

\begin{lemma}\label{nonbnmodel} The locus of non-BN general K3 surfaces in $\cK_{g}$ is a union of  the NL divisors $D_{d,n}^{g}$  satisfying 
$$\sqrt{2(g-1)n}< d \leq \min\{g-1, \frac{n+2}{2}+g-\frac{2g+2}{n+4}\}.$$
In particular, a quasi-polarized K3 surface $(S,L)$ of genus $6\leq g\leq 12 $ and $g\neq 11$  is non-BN general if and only if it lies in one of the following NL divisors:
\begin{enumerate}[$(I)$.]
\item $D^{g}_{d,0}$, $d=1,2,\ldots, [\frac{g-1}{2}]+1$;
\item  \begin{itemize}[-]
\item $g=6$, $ D^{6}_{5,2}$
\item $g=7$, $ D^{7}_{5,2}$, $D^{7}_{6,2}$
\item $g=8$, $ D^{8}_{6,2}$, $D^{8}_{7,2}$
\item $g=9$, $ D^{9}_{6,2}$, $D^{9}_{7,2}$
\item $g=10$, $ D^{10}_{7,2}$, $D^{10}_{8,2}$, $D^{10}_{9,4} $
\item $g=12$, $ D^{12}_{7,2}$, $D^{12}_{8,2}$, $D^{12}_{9,2}$, $D^{12}_{10,4}$, $D^{12}_{11,4}$
\end{itemize}
\end{enumerate} 
\end{lemma}
\begin{proof} 
First,  suppose $(S,L)$ is not BN general.  Then there exist line bundles $M,N\in \Pic(X)$ satisfying $M+N=L$ and 
\begin{equation}\label{nonbn}
h^0(M)h^0(N)\geq g+1
\end{equation}
To compute $h^0(M)$ and $h^0(N)$, let us recall some results about the linear systems on K3 surfaces. Let $F$ be an effective divisor on $S$. Saint-Donat has shown that
\begin{enumerate}[(i)]
\item When $|F|$ has no fixed component, then
\begin{itemize}
\item  $F^2>0$, then $|F|$ is base point free and $h^1(F)=0$. The Riemann-Roch Theorem yields
 $$h^0(F)=\frac{F^2}{2}+2.$$
\item  $F^2=0$, then $F=kE$ for some base point free line bundle $E$ satisfying $E^2=0$ and $h^0(F)=k+1$.
\end{itemize}
\item When $F=F'+\Gamma$ with the fixed component $\Gamma$ and $F'$ has no fixed component, then $F'$ is either base point free or $F'=kE$ by (i).
For each connected reduced fixed component $\Gamma'$ of F, we have $(\Gamma')^2=-2$ and $F'\cdot \Gamma'=0$ or $1$ if $F'$ is base point free; $E\cdot \Gamma'=0$ or $1$ if $F'=kE$.  The converse is also true.

A connected component $\Gamma'$ of $\Gamma$ is said to be of Type I if it is reduced with $\Gamma'\cdot F'=1$ (resp. $\Gamma'\cdot E=1$) and of Type II otherwise. The first cohomology of $F$ vanish  if $F$ only has Type I fixed component and $F'$ is base point free.
\end{enumerate}

Coming back to the proof, we discuss two cases: 
\\

\noindent {\bf Case 1}. If $|L|$  contains a fixed component,  then $S$ is an elliptic K3 surface with a section  which lies in $D^{g}_{1,0}$ (See also $\S$1.7). The assertion holds.
 \\
 
\noindent {\bf Case 2}. If $|L|$ has no fixed component which forces  $M^2\geq 0, N^2\geq 0$, we claim that we can find divisors $\widetilde{M},\widetilde{N}=L-\widetilde{M}$  such that $h^1(\widetilde{M})=h^1(\widetilde{N})=0$ and $h^0(\widetilde{M})h^0(\widetilde{N})\geq h^0(L)$. Admitting this, we set $\widetilde{M}^2=n$ and $L\cdot \widetilde{M}=d$, then $h^0(\widetilde{M})=\frac{n}{2}+2$ and $h^0(\widetilde{N})=g+\frac{n}{2}+1-d$.  Without loss of generality, we assume $h^0(\widetilde{M})\leq h^0(\widetilde{N})$. Thus we obtain 
\begin{equation}\label{ineuq}
\begin{cases}
     \left(\frac{n}{2}+2\right)\left(g+\frac{n}{2}+1-d\right)\geq g+1;& \text{ } \\
    g+\frac{n}{2}+1-d \geq \frac{n}{2}+2; & \text{} \\
     d^2-n(2g-2)>0. & 
\end{cases}
\end{equation}
Here the last inequality in \eqref{ineuq} is just the Hodge index theorem.  Thus $(S,L)$ is contained in $\cC_{d,n}^{g}$, where $d,n $ satisfies \eqref{ineuq}.

Now we prove the claim. Note that  the inequality $h^0(M)h^0(L-M)\geq h^0(L)$ still holds if we replace $M$ by its base point free part. So we can assume that $M$ is base point free. 

Denote by  $\Gamma_{I}$ ($\Gamma_{II}$) the sum of all Type I (Type II) fixed components of $N$. Then we can write $N=N'+\Gamma_I+\Gamma_{II}$ and $N'$ has no fixed component.   Then the line bundles
$$\widetilde{M}:=\begin{cases}
  M+\Gamma_{II}  & \text{if $N'$ is base point free}, \\
     M+(k-1)E+\Gamma_{II} & \text{if $N'=kE$ and $E^2=0$}.
\end{cases}$$
and $\widetilde{N}=L-\widetilde{M}$ has zero first cohomology. This is because $L$ is nef, which implies $\widetilde{M}$ and $\widetilde{N}$ have no Type II fixed component by a simple intersection computation. Thus we have proven the claim. 

Moreover, our assertion follows easily from the fact that the union of $\cC_{d,n}^g$ is the same as the union of $D_{d,n}^g$ for $d,n$ satisfying the condition \eqref{ineuq}.

Conversely, if $(S,L)$ lies in one of the NL divisors listed in the statement, then $\Pic(S)$ contains an element $M$ such that $M\cdot L=d$ and $M^2=n$ satisfying the given condition. We have the inequality $$h^0(M)h^0(L-M)\geq \left(\frac{M^2}{2}+2\right) \left(\frac{(L-M)^2}{2} +2\right)\geq g+1.$$  
\end{proof}

\begin{remark}
A similar computation can be found in \cite{JK04} Chapter 5, where Johnsen and Knutsen classify all non-Clifford general K3 surfaces of genus $6\leq g\leq 10$.  
\end{remark}

\subsection{Projective model of non-BN general K3 surfaces} In this subsection, we describe the  projective models of non-BN general K3 surfaces, which will be used later. Indeed, there is a projective model for the general K3 surface in each of NL divisors of Lemma 1.6.  Let us start with Saint-Donat's  result:
\begin{proposition}{\cite{Sa74}} If $\psi_L$ is not a birational to its image,  then 
\begin{enumerate}
  \item  $\psi_L$ is a generically $2:1$ map and $\psi_L(S)$ is a smooth rational normal scroll of degree $g-1$,  or a cone over a rational normal curve of degree $g+1$.
  \item $|L|$ has  a fixed component $D$, which is a smooth rational curve. The image of $\psi_{|L-D|}(S)$ is a rational normal curve of degree $g$ in $\PP^{g}$.
\end{enumerate}
\end{proposition}
The first (resp.~second) case happens if and only if $(S,L)$ admits an elliptic fibration with a section (resp.~bisection), which means that  $S$ contains an elliptic curve of degree one (resp.~two). Moreover, we have

\begin{lemma} 
For $2\leq g\leq 12$ and $g\neq 11$,  the K3 surface $(S,L)$ lies in $\cC^g_{2,0}$ if and only if  $\psi_L:S\rightarrow \PP^{g}$ is not birational to its image. Moreover, $\cC^g_{2,0}=D^{g}_{1,0}\cup D^{g}_{2,0}$ when $g\neq 7$, and $\cC^7_{2,0}=D^{g}_{1,0}\cup D^{g}_{2,0}\cup D^{7}_{5,2}$.
\end{lemma}

\begin{proof} By the discussion above,  $\psi_L$ is not birational if and only if it contains an elliptic curve of degree one or two.
The locus of those $(S,L)\in \cK_{g}$ containing an elliptic curve of degree one is the irreducible divisor $D^{g}_{1,1}$ because the corresponding rank two sublattice 
\begin{center}$\left(\begin{array}{cc}g& 1 \\1 & 0\end{array}\right)$
\end{center}
is always primitive. 

Similarly, for $6\leq g\leq 12$ and $g\neq 11$, the locus of those $(S,L)\in \cK_{g}$ containing an elliptic curve of degree two is the irreducible divisor $D^{g}_{2,1}$, except in the case $g=12$. When $g=12$, this locus is the union of the two irreducible divisors $D^{7}_{2,1}$ and $D^{7}_{5,2}$, 
because there are two primitive lattices 
\begin{center}$\left(\begin{array}{cc}12& 2 \\2 & 0\end{array}\right)$,~~~~and ~$\left(\begin{array}{cc}12&  5\\5 & 2\end{array}\right)$
\end{center}
which contain an element $\beta$ satisfying $\beta^2=0$ and $\beta\cdot L=2$.
\end{proof}

For the case $(S,L)$ is non-BN general and $\psi_L$ is birational, the image $\psi_L(S)$ is contained in some {\it rational normal scrolls} $\cT$, i.e. projective bundles over $\PP^1$. Here we say a rational number scroll $\cT$ is of type $(d_1,d_2,\ldots d_n)$ if 
$$\cT\cong\PP\left(\bigoplus\limits_{i=1}^n\cO_{\PP^1}(d_i)\right)$$
for $d_i\geq 0$.  The image of the natural morphism $\cT\rightarrow \PP \left(H^0(\bigoplus\limits_{i=1}^n\cO_{\PP^1}(d_i))\right)$ may be singular if $d_i=0$ for some $i$, and we say that $\cT$ is a singular rational normal scroll.

For the low genera appearing in this paper, Johnsen, Knutsen \cite{JK04} and Hana \cite{Ha02} have classified all the non-BN general K3 surfaces whose projective models lie in various rational normal scrolls. Using their results, one can find the projective models of the K3 surfaces lying in each irreducible component of the non-BN general locus.  
For instance, when $g=6$,  we have
\begin{itemize}
 \item  $(S,L)\in D^{6}_{3,0}$ if  $\psi_L(S)$  is the hypersurface of a rational normal scroll  of type $(2,1,1)$;
 \item   $(S,L)\in D^{6}_{5,2}$ if $\psi_L(S)$ is the complete intersection of a singular rational normal scroll  of type $(2,1,0,0)$. 
\end{itemize}
We refer the readers to \cite{JK04} Chapter $11$ and \cite{Ha02} Chapter $2$ for the complete list of all possible projective models.

\section{Birational models of $\cK_{g}$}
This section is devoted to the description of a birational model of $\cK_{g}$ via geometric invariant theory (GIT). More precisely, we can interpret the moduli space of BN general K3 surfaces as a moduli space of smooth complete intersections, using Theorem \ref{mukaimodel}. The latter can be constructed via GIT and this allows us to compute the Picard group of $\cK_{g}$ directly.

\subsection{GIT construction}   Recall that the general K3 surface in $\cK_{g}$ among the range of Mukai models is a complete intersection with respect to some vector bundle $E_{g}$ on a homogenous space $X_{g}$, by Theorem \ref{mukaimodel}. We obtain the GIT model of the moduli space as follows:

(I). For $6\leq g \leq 10$, let $V_{g}=H^0(X_{g}, \cO_{g}(1))$, the complete intersections $\psi_L(S)$ are parametrized by Grassmannians with a natural group action of $G_{g}$ coming from the action on $X_{g}$: 
\begin{itemize}
\item $W_{6}=\PP(H^0(F_{5},\cO_{\cF_5}(2) ))$, $G_{6}=PSL(2)$;
\item $W_{7}=Gr(8, V_{7})$,  $G_{7}= Spin(10)$;
\item $W_{8}=Gr(6, V_{8})$, $G_{8}=SL(6)$;
\item $W_{9}=Gr(4, V_{9})$,  $G_{9}=Sp(6)$;
\item $W_{10}=Gr(3, V_{10})$,  $G_{10}$ is the quotient of $G_2$ by its center;
\end{itemize} 
We denote by  $$\Phi_{g}:W_{g}^{ss}\q G_{g}\dashrightarrow \cK_{g}$$ the natural rational map to $\cK_{g}$, where $W_{g}^{ss}$ is the semistable locus of $W_{g}$. As proved by Mukai (cf.~\cite{Mu88} Theorem 0.2), $\Phi_{g}$ is birational, and thus the image of $\Phi_{g}$ contains an open subset of $\cK_{g}$. Our goal of this section is to describe the image of smooth complete intersections in $\cK_{g}$ via $\Phi_{g}$. 
 
(II).  The case of $g=12$ is slightly different. The BN general K3 surface is a hyperplane section of the smooth Fano threefold $Gr(3,V,N)\subseteq \PP^{13}$ for some non-degenerate $N$. This is a complete intersection with respect to the vector bundle $(\wedge^2\cF)^{\oplus 3}$ on $X_{12}=Gr(3,V)$. Since $H^0(Gr(3,V), \wedge^2\cF)\cong\wedge^2V^\vee$,  we have a natural parameter space $\cV_{12}$ of non-degenerate $Gr(3,V,N)$ which is birational to the GIT quotient $$Gr(3, \wedge^2V^\vee)^{ss}\q PGL(V).$$  The moduli space 
$\cK_{12}$ is birational to  the $\PP^{13}$-bundle $\cP_{12}\rightarrow \cV_{12}$,  where the fiber over an element $F\in\cV_{12}$ is the projective space $\PP(H^0(\cO_{F}(1)))$.
 In the rest of this section, we will discuss the GIT stability of the smooth elements in $W_g$. We start with the study of discriminant loci of complete intersections, which plays an important role in the proof of the main theorem.
 
\subsection{Discriminant loci of complete intersections}
First, we discuss general discriminant loci in moduli of complete intersections with respect to vector bundles, and then we restrict to our cases. 

Let $X$ be a smooth projective variety and $L$ a globally generated line bundle on $X$.  The {\it discriminant locus} $\Delta_L$ of $L$ is defined to the subset of $|L|$ parameterizing singular elements in $|L|$.  When $L$ is very ample,  we set $V_L=H^0(X,L)$ and let $(\PP^n)^\vee\cong \PP(V_L)$ be the dual projective space. The discriminant locus $\Delta_L$ is isomorphic to the dual variety $X^\vee$, which is an irreducible subvariety in $\PP(V_L)$. Moreover, one can easily get:
\begin{lemma}\label{dual}
The dual variety $X^\vee\cong \Delta_L$ is an irreducible hypersurface of $\PP(V_L)$ if there is an element $s'\in \Delta_L$ such that $(s')_0$ has only isolated singularities. 
\end{lemma}
\begin{proof}The proof is very standard, so we give only a sketch. Consider the incidence variety
\begin{equation}\label{incident}
\Sigma_L = \{(s,x)|(s)_0~ is~ singular~ at~ x \}\subseteq \PP(V_L)\times X,
\end{equation} with two projections 
$\pi_1:\Sigma_L\rightarrow \PP(V_L)$ and $\pi_2:\Sigma_L\rightarrow X$.

By our assumption, the image of the first projection $\pi_1$ is $\Sigma_L\neq \emptyset$ and $\pi_1$ is generically finite. Also, the second projection $\pi_2$ is surjective and each fiber of $\pi_2$ is a projective space $\PP^n$, where $n=\dim \PP(V_L)-\dim X-1$. Thus $\Delta_L=\pi_1(\Sigma_L)$ is irreducible of codimension one.
\end{proof}
\begin{remark} The irreducibility of the dual variety $X^\vee$ actually holds for any irreducible variety $X$. It is usually of codimension one in $\PP(V_L)$. When $X^\vee$ has codimension $r+1$ for some $r>0$, then $X$ is $\PP^r$-ruled (cf.~\cite{GKZ94} Chapter 1.1). 
\end{remark}

If $L$ is not very ample, we denote by $\varphi_L:X\rightarrow \PP(V_L)$ the morphism given by the linear system $|L|$ and let $Y$ be the image of $\varphi_L$.   The pullback induces an isomorphism 
\begin{equation}
\varphi_L^\ast: \PP(H^0(Y, \cO_Y(1)))\rightarrow \PP(V_L).
\end{equation}
Let $Z\subseteq Y$ be the union of the singular locus of the morphism $\varphi_L:X\rightarrow Y$ and the singular locus of $Y$. Suppose that there exists an element $s_0\in |\cO_Y(1)|$ whose zero locus has only isolated singularities lying outside $Z$. Then we have 
\begin{enumerate}
\item the discriminant locus $\Delta_{\cO_Y(1)}$ of $\cO_Y(1)$ is irreducible of codimension one.
\item the discriminant locus $\Delta_{L}$ contains an open subset of the image $\varphi_L^\ast (\Delta_{\cO_Y(1)})$, and hence contains $\varphi_L^\ast (\Delta_{\cO_Y(1)})$ as an irreducible component.
\end{enumerate}
It follows that 
\begin{proposition}
Suppose that there exists an element $s_0\in |\cO_Y(1)|$ whose zero locus has only isolated singularities lying outside $Z$.  There is an irreducible component  $\Delta^\circ_L\subseteq \Delta_L$, which is  isomorphic to  $\Delta_{\cO_Y(1)}$ and has codimension one in $\PP(V_L)$.  In particular, $\Delta_L^\circ=\Delta_L$ if the morphism $\varphi_L$ is smooth.
\end{proposition}

\begin{remark}
Without the assumption on the singularities of $(s_0)_0$, one can still obtain an irreducible subvariety of $\Delta_L$ isomorphic to $\Delta_{\cO(Y)(1)}$, but it may not be an irreducible component of $\Delta_L$. See Example 1.1.0 in \cite{LM08}. 
\end{remark}

Now, we extend the results above to the case of a globally generated vector bundle $E$ on $X$.  This will allow us to deal with  complete intersections. 
Let us  define
\begin{equation}
\begin{aligned}
\Delta_E:=\{s\in \PP( H^0(X,E))^\vee &|~(s)_0~\hbox{is either singular, or}  \\ & \hbox{not a complete intersection w.r.t. $E$}\},
\end{aligned}
\end{equation}
Then we have: 
 \begin{proposition}\label{discriminantvb}Let $\PP(E)$ be the associated projective bundle on $X$, and $\cO_{\PP(E)}(1)$ its relatively ample line bundle.  There is a natural isomorphism 
$$\Psi: H^0(\PP(E),\cO(1))\xrightarrow{\sim} H^0(X,E) .$$
The zero locus of a section $s \in H^0(\PP(E),\cO(1))$ in $\PP(E)$ is smooth if and only if the section $\Psi(s)\in H^0(X,E)$ is non-degenerate and the zero locus of $\Psi(s)$ is smooth in $\PP(E)$. In other words,
$\Delta_E$ is isomorphic to $\Delta_{\cO_{\PP(E)}(1)}$. 

Moreover, with the same assumption as in Proposition 2.4 for the line bundle $\cO_{\PP(E)}(1)$,  $\Delta_E$ contains an irreducible hypersurface $\Delta_E^\circ$ in $\PP(H^0(X,E))$. 
\end{proposition}
\begin{proof}
See \cite{Mu92} Proposition 1.9 for the first assertion. The last assertion follows directly from Proposition 2.4.
\end{proof}
If the vector bundle $E=E_0^{\oplus k} $ is a direct sum of $k$-copies of a globally generated vector bundle $E_0$ on $X$, we parametrize the complete intersections with respect to $E$ by the Grassmannian $Gr(k, V_{E_0})$ and define
\begin{equation}\label{discrimant2}
\begin{aligned}
\Delta_{k,E_0}=\{ & \left<s_1,\ldots, s_k\right>\in Gr(k, V_{E_0}) | (\oplus_i s_i)_0 \hbox{ is either singular} \\ & \hbox{or not a complete intersection}\},
\end{aligned}
\end{equation}
which is called the {\it discriminant locus} of $E_0^{\oplus k}$ in $Gr(k, V_{E_0})$.
\begin{lemma}\label{discriminant locus in gm}
The same assertion of Proposition 2.5 holds for $\Delta_{k,E_0}$ in $Gr(k,V_{E_0})$. 
\end{lemma}
\begin{proof}
One can use a similar argument as in Lemma \ref{dual}, where the fiber of the second projection will be an irreducible Schubert variety, or just note that $\Delta_{k,E_0}$ is the quotient $(\Delta_{E_0}\cap \PP(H^0(X,E))^s)/ SL(k)$. We omit the details here.
\end{proof} 

\begin{remark}
For our purpose, we actually only deal with the case where the vector bundle $E_0$ is very ample, so the discriminant locus $\Delta_{k,E_0}$ will be irreducible.  One may observe that the dual of universal  subbundle $\cF$ on $Gr(3,V)$ is not ample, but we have another construction for that case.  Proposition 2.5 and Lemma 2.6 are motivated by the study of K3 surfaces of genus $18$ and $20$ (cf. \cite{Mu92}).   
\end{remark}

\subsection{GIT Stability} Now we discuss GIT stability of the smooth complete intersections in $\S$2.1.  We say that a non-degenerate threefold $Gr(3,V,N)$ is {\it non-special} if it does not belong to Prokhorov's class of genus 12 Fano threefolds defined in \cite{Pr90}. Then we start with a useful lemma:
\begin{lemma}\label{finiteauto}
Let $X$ be either a K3 surface or a non-special smooth non-degenerate threefold $Gr(3,V,N)$. Let $L$ be the natural polarization on $X$. The group $Aut_L(X)$ of automorphisms $f:X\xrightarrow{\sim} X$ that preserve the polarization $L$, i.e.  $f^\ast L\cong L$, is finite.
\end{lemma}
\begin{proof}
See \cite{Hu} Chapter 5, Proposition 3.3 for K3 surfaces and \cite{Pr90} for automorphism groups of Fano threefolds. 
\end{proof}
\begin{remark}
By \cite{Pr90}, there are only three types of Fano threefolds of genus $12$ with infinite automorphism group. One is the Mukai-Umemura manifold constructed in \cite{MU83} with automorphism group $SL(2)$, and the other two have automorphism group $G_a$ and $G_m$. The moduli space for each of these types is at most one dimensional. 
\end{remark}

\begin{thm}\label{smstable}
For $6\leq g\leq 10$, if $S\in \cW_{g}$ is one of the smooth complete intersections described above, then $S$ is GIT stable in $\cW_g$.  Moreover, the smooth non-degenerate complete intersection $Gr(3,V,N)\subseteq Gr(3,V)$ is GIT stable in $Gr(3,\wedge^2 V^\vee)$ if it is non-special.
\end{thm}
\begin{proof}
Let $\Delta_{g}\subseteq W_{g}$ denote the discriminant locus, in the sense of \eqref{discrimant2}, for $6\leq g \leq 10$.  By Lemma \ref{discriminant locus in gm},  $\Delta_{g}$ is an irreducible hypersurface in $W_{g}$, since there are complete intersections with only one rational double point. Indeed, the K3 surface with a primitive Picard lattice \begin{center}$\left(\begin{array}{cc}2g-2& 0 \\0 & -2\end{array}\right)$
\end{center}
will be such a singular complete intersection.

The discriminant locus $\Delta_{g}$  is thus cut out by a single homogeneous equation  $\Omega_{g}=0$, named the {\it discriminant form}. Moreover, it is easy to see that the discriminant form $\Omega_{g}$ is $G_{g}$-invariant because the property of singularity is preserved under changes of coordinates. Hence $\Omega_{g}$ is  a $G_{g}$-invariant function. The semistability of the smooth surface $S$ then follows from the fact that $\Omega_{g}$ does not vanish at $S$. It follows that $S$ is GIT stable  because $S$ only has finite stabilizer by Lemma \ref{finiteauto}.

For the case $g=12$,  we employ a similar idea, using a discriminant form.  Let $\cW_{12}\subseteq Gr(3,\wedge^2 V^\vee)$ be the locus of non-degenerate elements in $Gr(3,\wedge^2 V)$ and  $\Delta:=Gr(3,\wedge^2V^\vee)\backslash \cW_{12}$  the complement of $\cW_{12}$.  As any smooth non-degnerate $Gr(3,V,N)$ has only finitely many automorphisms if it is non-special, it suffices to show that $\Delta$ is an irreducible $G_{12}$-invariant divisor in $Gr(3,\wedge^2V^\vee)$. 
To prove this, let us consider the incident variety:
$$\Omega=\{(N,[v])|\exists~ \omega\in N~s.t.~\omega\wedge v~is ~decomposable\}\subseteq Gr(3,\wedge^2V^\vee)\times \PP(V^\vee).$$
with the first projection $\pi_1(\Omega)=\Delta$.   The second projection $\pi_2:\Omega\rightarrow \PP(V^\vee)$ is surjective, and the fiber $\pi_2^{-1}(v)$ for any $v\in \PP(V^\vee)$ can be described as follows:  

Set $D_v=\{[\omega]\in \PP(\wedge^2 V^\vee)|~\omega\wedge v~is ~decomposable\}\subseteq \PP(\wedge^2 V^\vee) $, and $V_v=V^\vee/v$ the quotient vector space.  Then $D_v$ is isomorphic to the  $14$-dimensional  irreducible variety  $$(\widetilde{Gr}(2,  V_v)\times V_v)/\CC^\ast,$$
where  $\widetilde{Gr}(2,V_v)$ is the space of  decomposable vectors in $\wedge^2 V_v$, and $\CC^\ast $ acts simultaneously by scaling on $\widetilde{Gr}(2,V_v)$ and $V_v\cong v\wedge V_v$. 

Next, we define the correspondence 
 $$\xymatrix{ ~&\cR_v=\{(N, \omega)|~\omega\subset N\}\subseteq\ar[dl]_{p}\ar[dr]^q Gr(3,\wedge^2 V^\vee)\times D_v & ~\\ Gr(3,\wedge^2 V^\vee)  & & D_v}$$
where the fiber of the second projection $q$ over $\omega\in D_v $ is isomorphic to the Grassmannian $Gr(2, \wedge^2 V^\vee/\omega)$.  Now we claim that the first projection $p$ is generically finite,  and then we have $\pi_2^{-1}(v)=p(\cR_v)$ is irreducible of dimension $50$.

To prove the claim, we need only show the existence of a fiber of $p$ which is finite. Actually, it is easy to see that there exists a three dimensional subspace $N\subseteq \wedge^2 V_v^\vee$ such that the plane $\PP(N)$ meets $Gr(2,V_v^\vee)$ in $\PP(\wedge^2V_v^\vee)$ at finitely many points. Identifying $N$ as an element in $Gr(3,\wedge^2V^\vee)$, we know that the fiber $p^{-1}(N)$ is just the intersection of $\PP(N)$ and $Gr(2,V_v^\vee)$ in $\PP(\wedge^2V_v^\vee)$, which is finite. Thus, we have proven the claim, and we have that $\Omega$ is irreducible of dimension $56$.

For the rest of the assertion, we still need to understand the general fiber of the projection $\pi_1:\Omega\rightarrow Gr(3,\wedge^2 V^\vee)$.  Indeed, for general $N\in \Delta$, it contains only one 1-dimensional linear subspace spanned by an element $\omega$ of the form $$\omega=f_1\wedge f_2+v\wedge f_3,$$ for some $v\in V^\vee$ and $f_1,f_2,f_3\in V_v$ are linearly independent. Then the fiber $\pi_1^{-1}(N)$ is just the three dimensional projective space $\PP(\left<v,f_1,f_2,f_3\right>)$.

Combining everything, we get  $\Delta=\pi_1(\Omega)$ is irreducible of dimension $53$ in $Gr(3,\wedge^2 V^\vee)$, and this completes  the proof.
\end{proof}

\section{Proof of the main theorem}
According to Borcherds' theta lifting theory, a beautiful result of Noether-Lefschetz theory on K3 surfaces is that the generating series of NL divisors (in the sense of \cite{MP07}) is a vector-valued modular form (cf.~\cite{Bo99},\cite{MP07}).  As an application, Bruinier \cite{Br02} computes the dimension of the subspace $\Pic_\QQ^{NL}(\cK_g)$ of $\Pic_\QQ(\cK_g)$ generated by NL divisors (see also \cite{LT13}). 

\begin{lemma}[\cite{LT13}]\label{dimformula}
Let $\rho_{g}$ be the dimension  of $\Pic^{NL}_\QQ(\cK_{g})$. Then
\begin{equation}\rho_{g}= \frac{31g+24}{24}-\frac{1}{4}\alpha_g-\frac{1}{6}\beta_g -  \sum\limits_{k=0}^{g-1}\left\{\frac{k^2}{4g-4}\right\}-\sharp\left\{k~|~\frac{k^2}{4g-4}\in\ZZ, 0\leq k\leq g-1 \right\} 
\end{equation}
where 
$$\alpha_g= \begin{cases}
     0, & \text{ if $g$ is even }; \\
      \left( \frac{2g-2}{2\ell-3}\right)& \text{otherwise}.
\end{cases}, ~~~~\beta_g= \begin{cases}
     \left( \frac{g-1}{4g-5}\right)-1, & \text{ if $g\equiv 1\mod 3$}, \\
   \left( \frac{g-1}{4g-5}\right)+\left( \frac{g-1}{3}\right)   & \text{otherwise}.
\end{cases}$$ and $\left( \frac{a}{b}\right)$ is the Jacobi symbol. In particular, $\rho_{g}=6,7,7,8,9,12 $ for $g=6,7,8,9, 10,11$.
\end{lemma}

Now we are ready to give the proof of our main theorem:
\begin{proof}We still separate the proof into two cases: (I) $6\leq g\leq 10$ and (II) $g=12$. For case (I), we let $\cU_{g}\subseteq W_{g}$ be the open subset parametrizing the smooth complete intersections.  By definition, $\cU_{g}$ is the complement of the discriminant locus $\Delta_{g}$ in $W_{g}$. By Theorem \ref{smstable}, the reductive Lie group $G_{g}$ acts properly on $\cU_{g}$. Then  we get an open immersion: $$\Phi'_{g}:\cU_{g}/G_{g}\longrightarrow \cK_{g},$$ 
and the image $\Phi_{g}'(\cU_{g}/G_{g})$ is the complement of the union of NL divisors listed in Theorem 1.1.

Moreover, note that  the complement  $W_{g}\backslash \cU_{g}$ is an irreducible divisor $\Delta_{g}$, the Picard group of $\cU_{g}$ with rational coefficients has to be trivial. Let $\Pic_\QQ(\cU_{g})_{G_{g}}$ be the group of $G_{g}$-linearized line bundles on $\cU_{g}$. By \cite{KKV89} Proposition 4.2,  there is an injection
\begin{equation}
\Pic(\cU_{g}/G_{g})\hookrightarrow \Pic(\cU_{g})_{G_{g}}.
\end{equation}
The forgetful map  $\Pic_\QQ(\cU_{g})_{G_g}\rightarrow \Pic_\QQ(\cU_g)$ has kernel the group of rational characters $\chi(G_g)$, which is trivial because $G_g$ is simple. Thus $\Pic_\QQ(\cU_g/G_g)$ is trivial, and the Picard group $\Pic_\QQ(\cK_{g})$ is spanned by the boundary divisors $\cK_{g}\backslash {\rm Im}(\Phi_g')$, which are precisely the NL divisors stated in the theorem.  Finally, these NL divisors form a basis of $\Pic_\QQ(\cK_{g})$ by the dimension counts of Lemma \ref{dimformula}.
 
For case (II), we let $\cU_{12}\subseteq \cP_{12}$ be the locus parametrizing complete intersections of a smooth non-special threefold $Gr(3,V,N )\subseteq \PP^{13}$ and a hyperplane section in $\PP^{13}$.  By Theorem \ref{smstable}, we know that $\cU_{12}$ admits a fibration to the orbit space 
\begin{equation}\label{nonspecialspace}
\cW^{ns}_{12}/ SL(V).
\end{equation}
where $\cW^{ns}_{12}\subseteq Gr(3,\wedge^2 V^\vee)$ is the parameter space of non-special threefolds $Gr(3,V,N)$.  Then the open subset $$\cW^{ns}_{12}\hookrightarrow Gr(3,\wedge^2 V^\vee)\backslash \Delta$$ has Picard number zero, and so does the quotient $\cW^{ns}_{12}/ SL(V)$. The fiber of $\cU_{12}\rightarrow \cW^{ns}_{12}/ SL(V)$ over a threefold $F\in \cW^{ns}_{12}/ SL(V)$ is the complement of the discriminant $\Delta_{\cO_F(1)}$ in $\PP(\cO_{F}(1))$. Therefore, $\dim \Pic_\QQ(\cU_{12})=0$.
 
Similarly, one observes that the image of $\cU_{12}\hookrightarrow \cK_{12}$ is exactly the complement of the NL divisors listed in the main theorem statement, union with a high codimension ($\geq 7$) subvariety (see Remark 8) parametrizing the K3 surfaces in special threefolds.  This shows that the Picard group is generated by these divisors in the boundary of $\cU_{12}$ in $\cK_{12}$.  \end{proof}
By Lemma \ref{dimformula}, the rank of $\Pic_\QQ(\cK_{12})$ is $11$, hence the generators in \eqref{generator} is not linearly independent. They will satisfy a linear relation given in Remark \ref{relation}.
\subsection{Further discussion}

\subsubsection{Elliptic divisors} Among all the Noether-Lefschetz divisors, there is an interesting subcollection of divisors  $\{D^{2l}_{d,0},d\geq 1\}$, called {\it elliptic} divisors. The union of the elliptic divisors parametrizes those K3 surfaces admitting an elliptic fibration.  Such divisors appear among the generators of the Picard group $\Pic_\QQ(\cK_g)$. For instance, for $g\leq 5$,  the Picard group $\Pic_\QQ(\cK_{g})$ is spanned by elliptic divisors:
\begin{itemize}
\item $g\leq4$: $\{D^{g}_{d,0}, d=0,1,\ldots,g-1 \}$;
\item $g=5$: $\{D^{5}_{d,0}, d=0,1,\ldots,3 \}$.
\end{itemize}
A natural question is
\begin{question}
Are all the NL divisors supported on elliptic divisors? Equivalently,  is the subgroup $\Pic^{NL}_\QQ(\cK_{g})$ generated by elliptic divisors $\{D_{d,0}^{g},d\in\NN\}$ in $\Pic_\QQ(\cK_g)$?
\end{question}

This question can be rephrased in terms of the coefficients of vector-valued modular forms in the sense of Borcherds. For example,  Maulik \cite{Mau12} has shown that the Hodge line bundle on $\cK_{g}$ is supported on elliptic divisors. His proof relies on an estimate of the coefficients of a vector-valued cusp form (cf.~\cite{Mau12} Lemma 3.7).

More precisely, the space of NL divisors modulo the Hodge bundle is isomorphic to the space of vector-valued cusp forms.  Each NL divisor $D^g_{d,n}$ (in the sense of \cite{MP07}) corresponds to a vector-valued Poincar\'e series, computed in \cite{Bru02} $\S$1.2, whose coefficients are interpreted as generalized Kloosterman sums.  In principal, one can compute the coefficients of all modular forms using this basis. Here we did not do the computation of these modular forms here, but it would be interesting to check the low genus cases in a subsequent paper.
\begin{remark}\label{relation}
Using the relation between NL-divisors and the space of vector-valued cusp forms of weight $\frac{21}{2}$, Arie Peterson has computed that 
\begin{equation}
3D_{8,2}^{12}-D_{9,2}^{12}-4D_{10,4}^{12}+2D_{11,4}^{12}+8D_{4,0}^{12}-5D_{5,0}^{12}+D_{6,0}^{12}=0.
\end{equation}
This result will appear in his thesis. 
\end{remark}

\subsubsection{Arithmetic group cohomology} 
 As shown in \cite{CK98} $\S$6.3 (see also \cite{LT13} for the more general case), the second cohomology group $H^2(\cD/\Gamma_g, \QQ)$ has a pure Hodge structure and we have an isomorphism
\begin{equation}\label{cohomology}
\Pic_\QQ(\cD/\Gamma_g)\xrightarrow{\sim} H^2(\cD/\Gamma_g, \QQ)
\end{equation}
from the exponential exact sequence. The latter is also isomorphic to the arithmetic group cohomology $H^2(\Gamma_g, \QQ)$. It follows that 
\begin{corollary}
The second Betti number of the arithmetic group $\Gamma_{g}$ is $6,7,7,8,9,12$ for $g=6,7,\ldots,10,12$.
\end{corollary}
From \eqref{cohomology}, we know that the Noether-Lefschetz conjecture is equivalent to the cohomology group of $\cD/\Gamma_g$ being spanned by (classes of) locally Hermitian symmetric subvarieties. Recently, there have emerged more standard ways to study this problem on locally Hermitian symmetric varieties of orthogonal type, cf.~\cite{HH07},\cite{BMM11}. In a sequel to this paper, we will approach the conjecture from this direction.

\bibliographystyle {plain}
\bibliography{Picard}

\end{document}